\documentclass[12pt,leqno]{amsart}
\usepackage{amssymb}
\usepackage[hidelinks]{hyperref}

\calclayout
\makeatletter
\renewcommand{\@settitle}{\begin{center}\large\bfseries\@title\end{center}}
\renewcommand{\section}{\@startsection{section}{1}{\parindent}
  {1.2\baselineskip plus .3\baselineskip}{-.5em}{\normalfont\bfseries}}
\makeatother
\hypersetup{pdftitle={The exponent in the Erdos-Fuchs theorem for k-fold sums},
  pdfauthor={Xiao-Hui Yan}}
\numberwithin{equation}{section}
\newtheorem{theorem}{Theorem}[section]
\newtheorem{corollary}[theorem]{Corollary}
\theoremstyle{remark}
\newtheorem*{remark}{Remark}

\begin{document}

\title[The optimal exponent for $k$-fold sums]{The optimal exponent in the
Erd\H{o}s--Fuchs theorem for $k$-fold sums}
\author{Xiao-Hui Yan}
\address{School of Mathematics and Statistics, Anhui Normal University,
Wuhu 241002, P. R. China}
\email{yanxiaohui\_1992@163.com}
\date{}

\begin{abstract}
Let $k\ge 2$ be an integer. For an infinite set $A$ of nonnegative integers and a real number $x\ge0$, let $r_k(A,x)$ be the number of $k$-tuples $(a_1,\ldots,a_k)\in A^k$ satisfying $a_1+\cdots +a_k\le x$.
In this paper, we prove that
$r_k(A,x)=cx+o(x^{(1-1/k)^2})$ is impossible for any real number $c>0$. The exponent $(1-1/k)^2$ is best possible for $k\ge3$ by a recent result of Guo and Yan. Moreover, we obtain a mean-square estimate of $r_k(A,x)-cx$ and confirm a conjecture of Chen and Tang for every $k\ge3$.
\end{abstract}
\maketitle

\section{Introduction}

Let $\mathbb N$ be the set of all nonnegative integers and $k\ge2$ be an integer. For an infinite set $A\subseteq \mathbb N$ and a real number $x\ge0$, let
\[
 r_k(A,x)=\#\{(a_{1},\ldots,a_k):a_{1}+\cdots+a_{k}\le x,\ a_{1},\ldots,a_k\in A\}.
\]

In 1956, Erd\H{o}s and Fuchs \cite{EF} proved that
\[
 r_2(A,x)=cx+o(x^{1/4}(\log x)^{-1/2})
\]
is impossible for every
$c>0$. Jurkat (unpublished) and later Montgomery and Vaughan \cite{MV}
removed the logarithmic factor $(\log x)^{-1/2}$ in the error term.

In 2009, Tang \cite{Tang} extended the Erd\H{o}s--Fuchs theorem to $k$-fold
sums by proving that
\begin{equation}\label{eq:Tang09}
 r_k(A,x)=cx+o(x^{1/4})
\end{equation}
cannot hold for any set $A\subseteq\mathbb N$ and any constant $c>0$. Later, Chen and Tang \cite{CT11} gave a quantitative version of this result.
The exponent $1/4$ is thus the same as in the binary case, regardless
of the number of summands. This raises the question of how large an
exponent can replace $1/4$ for a given integer $k$. 

In this paper, we improve this exponent to $(1-1/k)^2$.
\begin{theorem}\label{thm:point}
Let $k\ge2$ be an integer and $c>0$ be a real number. There is no infinite set
$A\subseteq\mathbb N$ such that
\[
 r_k(A,x)=cx+o\bigl(x^{(1-1/k)^2}\bigr).
\]
\end{theorem}

The exponent $(1-1/k)^2$ in Theorem \ref{thm:point} is optimal. Indeed, Guo and Yan \cite{GuoYan} recently obtained an inverse Erd\H{o}s--Fuchs theorem by proving that if $k\ge3$, then there exist an infinite set $A\subseteq\mathbb N$ and a real number $c>0$ such that
\[
 r_k(A,x)=cx+O\bigl(x^{(1-1/k)^2}\bigr).
\]
This improves the earlier bound $O\bigl(x^{1-1/k-1/k^2}\sqrt{\log x}\bigr)$ of Dai and Pan \cite{DP}.
For $k=2$, the best known upper bound for the error term is due to
Yang and Zhao \cite{YZ}, who improved the bound $O(x^{1/4}\log x)$ of Ruzsa~\cite{Ruzsa} to $O(x^{1/4}\sqrt{\log x})$.

For other generalizations of the Erd\H{o}s--Fuchs theorem and its inverse problem, see, for example, \cite{Bateman,Hayashi,RS,Rue,Vaughan,YL21,YL}.

For two real numbers $c>0$ and $x\ge0$, write
\[
 P_k(A,c,x)=r_{k}(A,x)-cx.
\]
Chen and Tang \cite{CT12} proved that
\begin{equation}\label{eq:conjecture}
 \int_0^X P_k(A,c,x)^2\,dx\gg_{A,k,c}X^{3/2}
\end{equation}
for infinitely many positive integers $X$. In particular, \eqref{eq:conjecture} implies that \eqref{eq:Tang09} cannot hold. Later, Dai and Pan \cite{DP20} extended this result to sums of different sets.

Chen and Tang \cite{CT12} further conjectured that \eqref{eq:conjecture}
holds for all sufficiently large $X$.
Under the assumption that $r_k(A,x)=cx+o(x)$, they \cite{CT12} proved that
\begin{equation}\label{eq:ChenTangcondition}
 \int_0^X P_k(A,c,x)^2\,dx\gg_{A,k,c}\left(\frac{X}{\log X}\right)^{3/2}
\end{equation}
for all sufficiently large $X$.

In this paper, we obtain a weighted mean-square estimate for $P_k(A,c,x)$.

\begin{theorem}\label{thm:weighted}
For every integer $k\ge2$, every $c>0$, and every infinite set
$A\subseteq\mathbb N$, we have
\begin{equation}\label{eq:weighted}
 \int_0^\infty P_k(A,c,x)^2 e^{-x/N}\,dx
 \gg_{k,c} N^{1+2(1-1/k)^2}
\end{equation}
for all sufficiently large $N$.
\end{theorem}

For the unweighted mean square, we have the following result.

\begin{corollary}\label{cor:ordinary}
For every integer $k\ge2$, every $c>0$, and every infinite set
$A\subseteq\mathbb N$, we have
\begin{equation}\label{eq:ordinary}
 \int_0^X P_k(A,c,x)^2\,dx
 \gg_{k,c}\left(\frac{X}{\log X}\right)^{1+2(1-1/k)^2}
\end{equation}
for all sufficiently large $X$.
\end{corollary}

\begin{remark} Corollary~\ref{cor:ordinary} also holds for a sequence $A$ of nonnegative integers.
\end{remark}

Corollary~\ref{cor:ordinary} proves Chen and Tang's conjectured bound \eqref{eq:conjecture}
for every integer $k\ge 3$. For $k=2$, Corollary~\ref{cor:ordinary}  gives \eqref{eq:ChenTangcondition} without
the assumption $r_2(A,x)=cx+o(x)$. However, the conjectured lower bound
of order $X^{3/2}$ remains open in this case.

For any real number $0<\delta\le1$, let
\[
 G_\delta=\left\{x\ge1:
 \int_0^x P_k(A,c,t)^2\,dt
 \ge \delta x^{1+2(1-1/k)^2}\right\}.
\]
As an application of Theorem~\ref{thm:weighted}, we next estimate
the lower logarithmic density of $G_\delta$.

\begin{corollary}\label{thm:small-scales}
For every integer $k\ge2$, every $c>0$, and every infinite set
$A\subseteq\mathbb N$,  there exists a constant $C_{k,c}>0$ such that
\[
  \liminf_{Y\to\infty}
  \frac{1}{\log Y}\int_{G_\delta\cap[1,Y]}\frac{dx}{x}
  \ge 1-C_{k,c}\,\delta^{1/(1+2(1-1/k)^2)}
\]
for every real number $0<\delta\le 1$.
\end{corollary}

\section{Proofs}\label{sec:weighted}

Let $k\ge2$ be an integer and $c>0$ be a real number and $A\subseteq \mathbb N$ be an infinite set. For $|z|<1$, let
\[
 F(z)=\sum_{a\in A}z^a,\qquad
 D(z)=\sum_{n\ge0}P_k(A,c,n)z^n.
\]
Then
\[
 \frac{F(z)^k}{1-z}=\sum_{n\ge0}r_k(A,n)z^n=\sum_{n\ge0}P_k(A,c,n)z^n+\sum_{n\ge0}cnz^n.
\]
Since
\[
 \sum_{n\ge0}nz^n=\frac{z}{(1-z)^2},
\]
it follows that
\begin{equation}\label{eq:identity}
 F(z)^k=\frac{cz}{1-z}+(1-z)D(z).
\end{equation}
Let $N\ge2$ and write
\[
 W(N)=\sum_{n\ge0}P_k(A,c,n)^2e^{-n/N},
 \qquad
 \Delta(t)=F(e^{-1/N+it})-F(e^{-2/N+it}).
\]

\noindent\textbf{Claim 1.} For $N^{-1}\le T\le1$,
\begin{equation}\label{eq:short-arc}
 \int_{-T}^T|\Delta(t)|^2\,dt
 \ll_{k,c} N^{2/k-1}+T^{4-2/k}W(N).
\end{equation}
\begin{proof}
For $|t|\le\pi$ and $1\le u\le2$, let
\[
z_u=e^{-u/N+it},\qquad d=N^{-1}+|t|,
\]
and
\[
 M(t)=|D(z_1)|+\int_1^2\left|\partial_u D(z_u)\right|\,du.
\]
Noting that $N\ge2$ and $1\le u\le 2$, we have $0<u/N\le 1$ and so
\[
e^{-u/N}\asymp1,\qquad
1-e^{-u/N}=\int_0^{u/N}e^{-v}\,dv\asymp\frac{u}{N}\asymp\frac1N.
\]
It follows that
\begin{equation}\label{eq:z_u}
|z_u|=e^{-u/N}\asymp1
\end{equation}
and
\begin{equation}\label{eq:1-z_u}
 |1-z_u|=\sqrt{(1-e^{-u/N})^2+2e^{-u/N}(1-\cos t)}\asymp \sqrt{N^{-2}+t^2}
\asymp d.
\end{equation}
For $1\le u\le2$, since
\[
D(z_u)=D(z_1)+\int_1^u\partial_v D(z_v)\,dv,
\]
it follows that
\begin{equation}\label{eq:D(z_u)}
 |D(z_u)|\le|D(z_1)|+\int_1^2\left|\partial_v D(z_v)\right|\,dv=M(t).
\end{equation}

If $|1-z_u|^2|D(z_u)|\le c|z_u|/2$ for every
$u\in[1,2]$, then \eqref{eq:identity}, \eqref{eq:z_u} and \eqref{eq:1-z_u} give
\[
|F(z_u)|^k
\ge \frac{c|z_u|}{|1-z_u|}-|(1-z_u)D(z_u)|
\ge \frac{c|z_u|}{2|1-z_u|}
\gg_c d^{-1},
\]
that is
\begin{equation}\label{eq:F(z_u)}
|F(z_u)|\gg_c d^{-1/k}.
\end{equation}
Differentiating \eqref{eq:identity} with respect to $u$, and using $\partial_u z_u=-z_u/N$, we obtain
\[
 kF(z_u)^{k-1}\partial_u F(z_u)
 =-\frac{cz_u}{N(1-z_u)^2}+\frac{z_u}{N}D(z_u)
   +(1-z_u)\partial_u D(z_u).
\]
It follows from \eqref{eq:z_u}, \eqref{eq:1-z_u}, \eqref{eq:D(z_u)}  and \eqref{eq:F(z_u)} that
\[
|\partial_u F(z_u)|
\ll_{k,c} d^{1-1/k}
\left(N^{-1}d^{-2}+N^{-1}M(t)+d\,|\partial_uD(z_u)|\right).
\]
Hence, by the definition of $M(t)$ and $d=N^{-1}+|t|\ge N^{-1}$,
\begin{align*}
 |\Delta(t)|
 &=|F(z_1)-F(z_2)|
 =\left|\int_1^2\partial_u F(z_u)\,du\right|
 \le\int_1^2\left|\partial_u F(z_u)\right|\,du\\
 &\ll_{k,c} d^{1-1/k}\left(N^{-1}d^{-2}+N^{-1}M(t)
             +d\int_1^2|\partial_u D(z_u)|\,du\right)\\
 &\ll_{k,c} N^{-1}d^{-1-1/k}+d^{2-1/k}M(t).
\end{align*}

If $|1-z_u|^2|D(z_u)|>c|z_u|/2$ for some $u=u_0\in[1,2]$, \eqref{eq:z_u}, \eqref{eq:1-z_u} and \eqref{eq:D(z_u)} give
\begin{equation}\label{eq:dM}
dM(t)
\gg |1-z_{u_0}||D(z_{u_0})|
>\frac{c|z_{u_0}|}{2|1-z_{u_0}|}
\gg_c d^{-1}.
\end{equation}
Hence, for every $u\in[1,2]$, \eqref{eq:identity} gives
\[
|F(z_u)|^k
\le \frac{c|z_u|}{|1-z_u|}+|1-z_u|\,|D(z_u)|
\ll_c d^{-1}+dM(t)
\ll_c dM(t).
\]
Noting that $1/k-1<0$, it follows from \eqref{eq:dM} that
\[
|F(z_u)|\ll_c (dM(t))^{1/k}=(dM(t))^{1/k-1}dM(t)\ll_c d^{2-1/k}M(t)
\]
for every $u\in[1,2]$. Hence,
\[
|\Delta(t)|\le |F(z_1)|+|F(z_2)|\ll_c d^{2-1/k}M(t).
\]

Therefore, in all cases, we have
\begin{equation}\label{eq:Dt}
|\Delta(t)|\ll_{k,c} N^{-1}d^{-1-1/k}+d^{2-1/k}M(t).
\end{equation}
Note that $z_u^n=e^{-un/N}e^{int}$ and
\[
D(z_u)=\sum_{n\ge0}P_k(A,c,n)e^{-un/N}e^{int},
\]
and so
\[
\partial_u D(z_u)=\sum_{n\ge0}P_k(A,c,n)\left(-\frac{n}{N}\right)e^{-un/N}e^{int}.
\]
By Parseval's identity,
\[
\int_{-\pi}^{\pi}|D(z_1)|^2\,dt
 =2\pi\sum_{n\ge0}P_k(A,c,n)^2e^{-2n/N}
\]
and
\begin{align*}
\int_{-\pi}^{\pi}|\partial_uD(z_u)|^2\,dt
 &=2\pi\sum_{n\ge0}P_k(A,c,n)^2
     \left(\frac nN\right)^2e^{-2un/N}\\
 &\le2\pi\sum_{n\ge0}P_k(A,c,n)^2
     \left(\frac nN\right)^2e^{-2n/N}.
\end{align*}
It follows from the Cauchy--Schwarz inequality that
\begin{align}\label{eq:llW(N)}
\int_{-\pi}^{\pi}M(t)^2\,dt
&\le 2\int_{-\pi}^{\pi}|D(z_1)|^2\,dt
 +2\int_{-\pi}^{\pi}\left(\int_1^2|\partial_uD(z_u)|\,du\right)^2\,dt\notag\\
&\le 2\int_{-\pi}^{\pi}|D(z_1)|^2\,dt+2\int_1^2\int_{-\pi}^{\pi}|\partial_uD(z_u)|^2\,dt\,du\notag\\
&\ll \sum_{n\ge0}P_k(A,c,n)^2
\left(1+\left(\frac nN\right)^2\right)e^{-2n/N}\notag\\
&\ll \sum_{n\ge0}P_k(A,c,n)^2e^{-n/N}=W(N),
\end{align}
where the last inequality follows from the fact that $(1+x^2)e^{-x}$ is
bounded for $x\ge 0$.
With $s=Nt$, we have
\begin{align}\label{eq:part1}
\int_{-T}^{T}N^{-2}d^{-2-2/k}\,dt
&=2N^{-2}\int_0^T(N^{-1}+t)^{-2-2/k}\,dt\notag\\
&=2N^{2/k-1}\int_0^{NT}(1+s)^{-2-2/k}\,ds\notag\\
&\le \frac{2}{1+2/k}\,N^{2/k-1}.
\end{align}
Noting that $N^{-1}\le T$, if $|t|\le T$, then
\[
 d=N^{-1}+|t|\le T+T=2T.
\]
It follows from $T\le1$, $4-2/k>0$ and \eqref{eq:llW(N)} that
\begin{equation}\label{eq:part2}
\int_{-T}^{T}d^{4-2/k}M(t)^2\,dt
\le (2T)^{4-2/k}\int_{-\pi}^{\pi}M(t)^2\,dt
\ll T^{4-2/k}W(N).
\end{equation}
It follows from \eqref{eq:Dt},\eqref{eq:part1} and \eqref{eq:part2} that
\begin{align*}
\int_{-T}^{T}|\Delta(t)|^2\,dt
&\ll_{k,c} \int_{-T}^{T}N^{-2}d^{-2-2/k}\,dt+\int_{-T}^{T}d^{4-2/k}M(t)^2\,dt\\
&\ll_{k,c}N^{2/k-1}+T^{4-2/k}W(N).
\end{align*}
\end{proof}

\noindent\textbf{Claim 2.}
For all sufficiently large $N$, we have
\begin{equation}\label{eq:W(N)lower}
 W(N)\gg_{k,c}N^{1+2(1-1/k)^2}.
\end{equation}
\begin{proof}
We may assume that
\begin{equation}\label{eq:W(N)}
W(N)\le N^{1+2(1-1/k)^2}.
\end{equation}
Note that
\[
 D(z)=\sum_{n\ge0}P_k(A,c,n)z^n,\qquad W(N)=\sum_{n\ge0}P_k(A,c,n)^2e^{-n/N}.
\]
For $j=2,3,4$, since
\[
\sum_{n\ge0}e^{-(2j-1)n/N}
=\frac{1}{1-e^{-(2j-1)/N}}
\asymp N,
\]
it follows from the Cauchy--Schwarz inequality and \eqref{eq:W(N)} that
\begin{align*}
 |D(e^{-j/N})|
 &=\left|\sum_{n\ge0}\left(P_k(A,c,n)e^{-n/(2N)}\right)e^{-(2j-1)n/(2N)}\right|\\
 &\le \left(\sum_{n\ge0}P_k(A,c,n)^2e^{-n/N}\right)^{1/2}
     \left(\sum_{n\ge0}e^{-(2j-1)n/N}\right)^{1/2}\\
 &= \frac{W(N)^{1/2}}{(1-e^{-(2j-1)/N})^{1/2}}
 \ll N^{1+(1-1/k)^2}.
\end{align*}
Since $0<1-e^{-j/N}\le j/N$, it follows that
\[
\bigl|(1-e^{-j/N})D(e^{-j/N})\bigr|
\ll N^{(1-1/k)^2}=o(N).
\]
Since $e^{j/N}-1=j/N+O(N^{-2})$, it follows that
\[
\frac{e^{-j/N}}{1-e^{-j/N}}
=\frac{1}{e^{j/N}-1}
=\frac Nj+O(1).
\]
Taking $z=e^{-j/N}$ in \eqref{eq:identity}, we have
\begin{align*}
F(e^{-j/N})^k
&=\frac{ce^{-j/N}}{1-e^{-j/N}}
 +(1-e^{-j/N})D(e^{-j/N})\\
&=\frac{cN}{j}+o(N)
=\frac{cN}{j}\bigl(1+o(1)\bigr).
\end{align*}
Then
\begin{align}\label{eq:b_a}
\sum_{a\in A}(e^{-a/N}-e^{-2a/N})^2
&=F(e^{-2/N})-2F(e^{-3/N})+F(e^{-4/N})\notag\\
&=(cN)^{1/k}
\left(2^{-1/k}-2\cdot3^{-1/k}+4^{-1/k}+o(1)\right)\notag\\
&\gg_{k,c} N^{1/k}.
\end{align}
Let $T>0$ and $h$ be real numbers. If $h\neq0$, then
\begin{align}\label{eq:kernel1}
\int_{-T}^{T}\left(1-\frac{|t|}{T}\right)e^{iht}\,dt
&=2\int_0^T\left(1-\frac{t}{T}\right)\cos(ht)\,dt\notag\\
&=\frac{2}{Th}\int_0^T\sin(ht)\,dt
=\frac{2(1-\cos(hT))}{Th^2}\notag\\
&=T\left(\frac{\sin(hT/2)}{hT/2}\right)^2\ge0.
\end{align}
If $h=0$, then
\begin{equation}\label{eq:kernel2}
\int_{-T}^{T}\left(1-\frac{|t|}{T}\right)e^{iht}\,dt
=2\int_0^T\left(1-\frac{t}{T}\right)\,dt=T>0.
\end{equation}
Write
\[
b_a=e^{-a/N}-e^{-2a/N}
=e^{-a/N}(1-e^{-a/N})\ge0.
\]
Then
\begin{align*}
\Delta(t)
&=F(e^{-1/N+it})-F(e^{-2/N+it})\\
&=\sum_{a\in A}e^{-a/N}e^{iat}
  -\sum_{a\in A}e^{-2a/N}e^{iat}\\
&=\sum_{a\in A}\bigl(e^{-a/N}-e^{-2a/N}\bigr)e^{iat}
=\sum_{a\in A}b_a e^{iat}
\end{align*}
and
\begin{align}\label{eq:Deltatlower}
\int_{-T}^{T}|\Delta(t)|^2\,dt
&\ge \int_{-T}^{T}\left(1-\frac{|t|}{T}\right)|\Delta(t)|^2\,dt\notag\\
&=\sum_{a,b\in A}b_ab_b
  \int_{-T}^{T}\left(1-\frac{|t|}{T}\right)e^{i(a-b)t}\,dt\notag\\
&\ge T\sum_{a\in A}b_a^2\gg_{k,c} TN^{1/k},
\end{align}
where the penultimate inequality follows from \eqref{eq:kernel1} and \eqref{eq:kernel2}, and the last inequality follows from \eqref{eq:b_a}.

For all sufficiently large $N$, let $T=B N^{1/k-1}$  such that $N^{-1}\le T\le1$, where $B$ is a sufficiently large constant depending only on $k,c$.
By \eqref{eq:short-arc} and \eqref{eq:Deltatlower}, we have
\[
 T N^{1/k}\ll_{k,c} N^{2/k-1}+T^{4-2/k}W(N),
\]
and so
\[
 W(N)\gg_{k,c} N^{2/k-1+(1-1/k)(4-2/k)}=N^{1+2(1-1/k)^2}.
\]
\end{proof}
\begin{proof}[Proof of Theorem~\ref{thm:weighted}]
For every integer $n\ge0$ and $0\le u<1$, we have
\[
r_k(A,n+u)=r_k(A,n),
\]
and so
\[
P_k(A,c,n+u)=r_k(A,n+u)-c(n+u)=P_k(A,c,n)-cu.
\]
With $x=n+u$, we have
 \begin{align}\label{eq:cell-integral}
\int_n^{n+1}P_k(A,c,x)^2\,dx
&=\int_0^1\left(P_k(A,c,n)-cu\right)^2\,du\notag\\
&=P_k(A,c,n)^2-cP_k(A,c,n)+\frac{c^2}{3}\notag\\
&=\frac14 P_k(A,c,n)^2
    +\frac34\left(P_k(A,c,n)-\frac{2c}{3}\right)^2\notag\\
&\ge\frac14 P_k(A,c,n)^2.
\end{align}
Since $e^{-x/N}\ge e^{-(n+1)/N}$ for any $x\in[n,n+1]$, it follows from \eqref{eq:cell-integral} and \eqref{eq:W(N)lower} that
\begin{align*}
\int_0^\infty P_k(A,c,x)^2e^{-x/N}\,dx
&\ge \sum_{n\ge0}e^{-(n+1)/N}
       \int_n^{n+1}P_k(A,c,x)^2\,dx\\
&\ge \frac{e^{-1/N}}4
       \sum_{n\ge0}P_k(A,c,n)^2e^{-n/N}\\
&=\frac{e^{-1/N}}4W(N)
 \gg_{k,c}N^{1+2(1-1/k)^2}.
\end{align*}
This completes the proof of Theorem \ref{thm:weighted}.
\end{proof}

\begin{proof}[Proof of Theorem~\ref{thm:point}]
Suppose that there exists an infinite set $A\subseteq\mathbb N$ such that
\[
P_k(A,c,x)=o(x^{(1-1/k)^2}).
\]
Noting that
\[
0<\int_0^\infty u^{2(1-1/k)^2}e^{-u}\,du<\infty,
\]
with $x=Nu$, we have
\begin{align*}
 \int_0^\infty P_k(A,c,x)^2e^{-x/N}\,dx
 &=o\left(\int_{0}^\infty x^{2(1-1/k)^2}e^{-x/N}\,dx\right)\\
 &=o\left(N^{1+2(1-1/k)^2}
 \int_0^\infty u^{2(1-1/k)^2}e^{-u}\,du \right)\\
  &=o\left(N^{1+2(1-1/k)^2}\right),
\end{align*}
a contradiction with \eqref{eq:weighted}. This completes the proof of Theorem \ref{thm:point}.
\end{proof}

\begin{proof}[Proof of Corollary~\ref{cor:ordinary}]
Let $X$ be sufficiently large. If $r_k(A,X/2)>2cX$, then, for any $x\in[X/2,X]$,
\[
P_k(A,c,x)=r_k(A,x)-cx\ge r_k(A,X/2)-cX>2cX-cX=cX,
\]
and so
\[
\int_0^X P_k(A,c,x)^2\,dx
\ge \int_{X/2}^X P_k(A,c,x)^2\,dx
\ge \frac{c^2}{2}X^3
\gg_{k,c} \left(\frac{X}{\log X}\right)^{1+2(1-1/k)^2}.
\]

Suppose now that $r_k(A,X/2)\le2cX$. Let 
\[
Y=X/(2k),\quad
A_Y=A\cap[0,Y],\quad
B=A_Y\cup\{m^k:m\in \mathbb{N},\  m^k>Y\}.
\]
Then 
\[
 (\#A_Y)^k\le r_k(A,X/2)\le 2cX=4kcY,
\]
$B\subseteq \mathbb N$ is an infinite set and 
\[
P_k(B,c,t)=P_k(A,c,t)
\]
for all $0\le t\le Y$. For any $t>Y$, we have
\[
 \#(B\cap[0,t])
 \le \#A_Y+t^{1/k}
 \le (4kcY)^{1/k}+t^{1/k}
 \ll_{k,c}t^{1/k}.
\]
Hence 
\[
r_k(B,t)\le(\#(B\cap[0,t]))^k\ll_{k,c}t,
\]
and so
\begin{equation}\label{eq:prefix}
|P_k(B,c,t)|\le r_k(B,t)+ct \ll_{k,c}t
\end{equation}
for all $t>Y$.

Let $N=X/(16k\log(2X))$. Note that
\[
\frac{Y}{2N}=4\log(2X),
\quad e^{-Y/(2N)}=(2X)^{-4},\quad 0<\int_0^\infty u^2e^{-u}\,du<\infty.
\]
By \eqref{eq:prefix},
\begin{align*}
 \int_Y^\infty P_k(B,c,t)^2e^{-t/N}\,dt
 &\ll_{k,c} e^{-Y/(2N)}\int_0^\infty t^2 e^{-t/(2N)}\,dt\\
 &=(2X)^{-4}(2N)^3\int_0^\infty u^2 e^{-u}\,du\\
 &\ll_{k,c}X^{-1},
\end{align*}
where the equality follows by setting $t=2Nu$. It follows from \eqref{eq:weighted} that
\[
\begin{aligned}
\int_0^X P_k(A,c,t)^2\,dt
 &\ge \int_0^Y P_k(A,c,t)^2e^{-t/N}\,dt=\int_0^Y P_k(B,c,t)^2e^{-t/N}\,dt\\
&=\int_0^\infty P_k(B,c,t)^2e^{-t/N}\,dt
  -\int_Y^\infty P_k(B,c,t)^2e^{-t/N}\,dt\\
&\gg_{k,c} N^{1+2(1-1/k)^2}\gg_{k,c} \left(\frac{X}{\log X}\right)^{1+2(1-1/k)^2}.
\end{aligned}
\]
This completes the proof of Corollary~\ref{cor:ordinary}.
\end{proof}

\begin{proof}[Proof of Corollary~\ref{thm:small-scales}]
Let 
\[
S(x)=\int_0^x P_k(A,c,t)^2\,dt.
\]
Since
\[
 r_k(A,t)=\#\{(a_{1},\ldots,a_k):a_{1}+\cdots+a_{k}\le t,\ a_{1},\ldots,a_k\in A\}\le (1+t)^k,
\]
it follows that
\[
|P_k(A,c,t)|
\le r_k(A,t)+ct
\le (1+t)^k+ct
\le (1+c)(1+t)^k
\]
and
\[
 S(x)\le (1+c)^2\int_0^x (1+t)^{2k}\,dt \ll_{k,c}(1+x)^{2k+1}.
\]
For any integer $m\ge0$, let
\[
u_m=\log(1+S(2^m)).
\]
Then the sequence $(u_m)_{m\ge0}$ is nondecreasing and
\[
 u_m \ll_{k,c} m+1.
\]
Let $H\ge H_0(k,c)$ be a sufficiently large constant in terms of $k$ and $c$. For every integer $M\ge1$, let
\[
T(M)=\bigl\{0\le m\le M:
 u_{m+r}-u_m>Hr \text{ for some integer } r\ge1\bigr\}.
\]
If $0\le m\le M$ and $m+r>2M$, then $r>2M-m\ge m$ and 
\[
u_{m+r}-u_m
\le u_{m+r}
\ll_{k,c} m+r+1
\ll_{k,c} 2r
\ll_{k,c} r.
\]
Hence, if $m\in T(M)$, then every $r$ for which $u_{m+r}-u_m>Hr$ satisfies $m+r\le 2M$. 

Let $m_0$ be the smallest integer in $T(M)$ and $r_0\ge1$ be an integer such that $u_{m_0+r_0}-u_{m_0}>Hr_0$. Let $m_{i+1}$ be the smallest integer with
\[
m_{i+1}\in T(M)\setminus\left(\bigcup_{j=0}^{i}\left[m_j,m_j+r_j\right)\right)
\]
and $r_{i+1}\ge1$ be an integer such that 
\[
u_{m_{i+1}+r_{i+1}}-u_{m_{i+1}}>Hr_{i+1}.
\]
Then
\begin{equation}\label{eq:boundTM}
H\sum_i r_i
<
\sum_i\bigl(u_{m_i+r_i}-u_{m_i}\bigr)
\le u_{2M}-u_0
\ll_{k,c}M
\end{equation}
and
\[\#T(M)\le\sum_i r_i\ll_{k,c}\frac{M}{H}.\]

Let $0\le m \le M$ be an integer with $m\notin T(M)$ and $N=2^m/(3H)$. Then, for every integer $r\ge1$,
\[u_{m+r}-u_m=\log\frac{1+S(2^{m+r})}{1+S(2^m)}\le Hr\]
and so
\begin{equation}\label{eq:Sm+r}
S(2^{m+r})
\le \bigl(1+S(2^m)\bigr)e^{Hr}.
\end{equation}
Since
\[
\sum_{j\ge0}e^{H(j+1)-3H2^j}
\le\sum_{j\ge0}e^{-2H(j+1)}
=\frac{1}{e^{2H}-1}<1,
\]
it follows from \eqref{eq:Sm+r} with $r=j+1$ that
\begin{align*}
 \int_0^\infty P_k(A,c,t)^2e^{-t/N}\,dt
 &=\int_0^{2^m} P_k(A,c,t)^2e^{-t/N}\,dt+\sum_{j\ge0}\int_{2^{m+j}}^{2^{m+j+1}} P_k(A,c,t)^2e^{-t/N}\,dt\\
 &\le S(2^m)+\sum_{j\ge0}e^{-2^{m+j}/N}S(2^{m+j+1})\\
 &= S(2^m)+\sum_{j\ge0}e^{-3H2^j}S(2^{m+j+1})\\
 &\le (1+S(2^m))
       \left(1+\sum_{j\ge0}e^{H(j+1)-3H2^j}\right)\\
 &\le 2(1+S(2^m)).
\end{align*}
It follows from Theorem~\ref{thm:weighted} that
\[
2\bigl(1+S(2^m)\bigr)
\ge \int_0^\infty P_k(A,c,t)^2e^{-t/N}\,dt
\gg_{k,c}N^{1+2(1-1/k)^2},
\]
that is,
\[
S(2^m)\gg_{k,c}(2^m/H)^{1+2(1-1/k)^2}.
\]
Then, for $2^m\le x<2^{m+1}$, there exists a constant $b_{k,c}>0$ such that
\[
 S(x)\ge S(2^m) \ge b_{k,c}(x/H)^{1+2(1-1/k)^2}.
\]
Let
\[
\delta_0=\min\{1,b_{k,c}H_0^{-(1+2(1-1/k)^2)}\}.
\]
For $0<\delta<\delta_0$, take
\[
H=(b_{k,c}/\delta)^{1/(1+2(1-1/k)^2)}.
\]
Then $H\ge H_0$ and $b_{k,c}H^{-(1+2(1-1/k)^2)}=\delta$. Thus there exists an integer $m_*=m_*(k,c,\delta)\ge 0$ such that,
whenever $m_*\le m\le M$ and $m\notin T(M)$, we have
\[
S(x)\ge \delta x^{1+2(1-1/k)^2}
\qquad (2^m\le x<2^{m+1}).
\]
In particular, $[2^m,2^{m+1})\subseteq G_\delta$ for every such $m$.

For $Y\ge2$, taking $M=\lfloor\log_2 Y\rfloor$ and using \eqref{eq:boundTM},
we obtain
\[
\int_{[1,Y]\setminus G_\delta}\frac{dx}{x}
\le \bigl(m_*+\#T(M)\bigr)\log 2
\le m_*\log 2+D_{k,c}\frac{\log Y}{H},
\]
where $D_{k,c}>0$ depends only on $k$ and $c$.
Hence,
\begin{align*}
\liminf_{Y\to\infty}
\frac{1}{\log Y}\int_{G_\delta\cap[1,Y]}\frac{dx}{x}
&=1-\limsup_{Y\to\infty}
\frac{1}{\log Y}\int_{[1,Y]\setminus G_\delta}\frac{dx}{x}\\
&\ge 1-\frac{D_{k,c}}{H}
 =1-D_{k,c}b_{k,c}^{-1/(1+2(1-1/k)^2)}\delta^{1/(1+2(1-1/k)^2)}.
\end{align*}

Finally, take
\[
C_{k,c}=\max\{D_{k,c}b_{k,c}^{-1/(1+2(1-1/k)^2)},\,\delta_0^{-1/(1+2(1-1/k)^2)}\}.
\]
The preceding estimate gives the desired bound for $0<\delta<\delta_0$.
For $\delta_0\le\delta\le 1$, we have
$1-C_{k,c}\delta^{1/(1+2(1-1/k)^2)}\le 0$, so the bound follows from the
nonnegativity of the lower logarithmic density.
This completes the proof of Corollary~\ref{thm:small-scales}.
\end{proof}

\medskip
\noindent\textbf{Acknowledgement.} This work was supported by the National Natural Science Foundation of China (Grant No.~12371005).

\end{document}